\documentclass[12pt]{article}
\usepackage{color}

\usepackage{amsmath}
\usepackage{amsfonts}
\usepackage{amsthm}
\usepackage{amssymb}
\usepackage{bm}
\usepackage[applemac]{inputenc}
\usepackage{enumerate}
\usepackage{stmaryrd}
\makeatletter

\newtheoremstyle{mystyle}{}{}{\slshape}{2pt}{\scshape}{.}{ }{} 

\newtheorem{thm}{Theorem}[section]
\newtheorem{cor}[thm]{Corollary}

\newtheorem{prop}[thm]{Proposition}
\newtheorem{lemme}[thm]{Lemma}

\newtheorem{conclusion}[thm]{Conclusion}

\theoremstyle{definition}

\theoremstyle{mystyle}

\theoremstyle{remark}
\newtheorem{rem}[thm]{Remark}

\DeclareMathOperator{\acl}{acl}

\title{Adding linear orders}
\author{S. Shelah \footnote{The author would like to thank the Israel Science Foundation for
 partial support of this research (Grant no. 710/07).
 Publication 979 on Shelah's list.} , P. Simon}

\begin{document}
\maketitle

\begin{abstract}We address the following question: Can we expand an NIP theory by adding a linear order such that the expansion is still NIP? Easily, if acl($A$)=$A$ for all $A$, then this is true. Otherwise, we give counterexamples. More precisely, there is a totally categorical theory for which every expansion by a linear order has IP. There is also an $\omega$-stable NDOP theory for which every expansion by a linear order interprets bounded arithmetic.
\end{abstract}

A well known open question is whether every unstable $NIP$ theory interprets an infinite linear order. We are concerned here with a question somewhat in the same spirit but going in a different direction: Can we expand an $NIP$ theory by adding a linear order on the whole universe so that the resulting theory is still $NIP$? We give a negative answer in two strong forms: 

1) There is an $\omega$-stable $NDOP$ theory of depth 2 for which every expansion by a linear order interprets bounded arithmetic (see section \ref{bbarith}).

2) There is a totally categorical theory for which every expansion by a linear order has $IP$.
\\

In the first section, we mention a few positive statements that are true (and easy): if $M$ is $NIP$ and $\acl(A)=A$ for all $A\subset M$, then $M$ can be linearly ordered so as to stay $NIP$. Also if $M$ is $\omega$-categorical, then $M$ can be linearly ordered so as to stay $\omega$-categorical (a well known fact) so we cannot expect to get the strong conclusion of 1) with an $\omega$-categorical theory.

Let us also note that, as proved in \cite{CS} (using results from \cite{BB}), adding a predicate for a small dense indiscernible sequence preserves $NIP$. In particular any unstable $NIP$ theory $T$ has an $NIP$ expansion which defines an infinite linear order.
\\

As far as we know, the question we address was first asked by Artem Chernikov. It came up again in discussions with Udi Hrushovski, which led to this work. We would like to thank him for that and for helping the second author understand the results of the first author.

\section{The easy case}

We fix a one-sorted relational language $\mathcal L$ not containing the symbol $<$, let $\mathcal L_<$ be the language with a unique binary relation $<$ and let $\mathcal L'=\mathcal L\cup \mathcal L_<$. Let $T$ be a complete $\mathcal L$-theory that eliminates quantifiers. Let $T'$ be the $\mathcal L'$-theory generated by $T$ and axioms stating that $<$ defines a dense linear order with no end-points.

We show in this section that if $T$ eliminates $\exists^\infty$, then $T'$ has a model-companion. Apart maybe from Proposition \ref{expand}, everything here is well known. We follow the exposition of \cite{TsuboiAmalg} which contains exactly what we need. The $\omega$-categorical case was already observed by Schmerl in \cite{Schmerl}.

Let $M$ be any structure, and $A\subset M$ a finite subset. We say that a formula $\phi(x_1,...,x_n)$ with parameters in $A$ is \emph{large} if it has a solution $a_1,...,a_n$ such that for all $i\neq j$, $a_i \notin \acl(A)$ and $a_i \neq a_j$. If $\phi(x_1,...,x_n;y)$ is a formula, then using elimination of $\exists^\infty$, it can be checked that the set of $d$'s such that $\phi(x_1,...,x_n;d)$ is large is definable (Fact 2 of \cite{TsuboiAmalg}).

\begin{prop}
Let $T$ be any $\mathcal L$-theory that eliminates quantifiers and $\exists^\infty$. Then $T'$ admits a model-companion $\tilde T$ axiomatized by $T$ along with the statements saying that for every large $\mathcal L$-formula $\phi(x_1,..,x_n)$ and quantifier-free large $\mathcal L_<$-formula $\theta(x_1,...,x_n)$, the conjunction $\phi(x_1,...,x_n)\wedge \theta(x_1,...,x_n)$ has a solution.

Furthermore, in $\tilde T$, the type of an algebraically closed set is given by the union of its $\mathcal L$-type and its $\mathcal L_<$-quantifier-free-type.
\end{prop}
\begin{proof}
The first statement is a special case of Lemma 8 of \cite{TsuboiAmalg}. The proof is straightforward: any model of $T'$ embeds into a model of $\tilde T$ by iteratively adding solutions to formulas. Existential closeness is clear.

The second statement is by easy back-and-forth.
\end{proof}

\begin{cor}
Let $M\models T$ be an $\mathcal L$-structure, and assume that algebraic closure is trivial in $M$: $\acl(A)=A$ for all $A\subset M$. Assume that $M$ admits elimination of quantifiers in $\mathcal L$. Then there is an expansion $T'$ of $T$ to $\mathcal L'=\mathcal L\cup \{<\}$ such that $<$ defines a linear order and $T'$ has elimination of quantifiers in $\mathcal L'$.
\end{cor}
\begin{proof}
Take $\tilde T$ as above. As $\acl(A)=A$ for all $A\subset M$, the type of any set $A$ is given by its $\mathcal L$-type and its quantifier-free $\mathcal L_<$-type. As $T$ eliminates quantifiers in $\mathcal L$, $\tilde T$ eliminates quantifiers in $\mathcal L'$.
\end{proof}

\begin{conclusion}
\begin{enumerate}
\item If $T$ is $\omega$-categorical, then $T$ admits an $\omega$-categorical expansion to $\mathcal L'$ in which $<$ defines a linear order.
\item If $T$ has trivial algebraic closure and is $NIP$, then it admits an $NIP$ expansion to $\mathcal L'$ in which $<$ defines a linear order.
\end{enumerate}
\end{conclusion}
\begin{proof}
The second point follows from the Corollary.

For the first point, notice first that any $\omega$-categorical theory eliminates $\exists^\infty$. More precisely, for finite $A$, $\acl(A)$ is finite. So for a given integer $n$, the type of an $n$-tuple $\bar a=(a_1,...,a_n)$ in a model of $\tilde T$ is given by its $\mathcal L$-type and the ordering of $\acl(\bar a)$. We see that there are only finitely many possibilities. Therefore $\tilde T$ is $\omega$-categorical.
\end{proof}

We end this section with a small observation.

\begin{prop}\label{expand}
Assume $T$ eliminates $\exists^\infty$. Let $M\models T$. If all infinite definable sets of $M$ has the same cardinality $\lambda$, then $M$ admits an expansion to a model $M'$ of $\tilde T$ as defined above.
\end{prop}
\begin{proof}
Fix a $\lambda$ saturated dense linear order $(\Omega,<_\Omega)$ without end points. We will build an injection $f:M\rightarrow \Omega$. Fix an enumeration $(\bar a_\alpha : \alpha<\lambda)$ of all pairs $\bar a_\alpha = (\phi_\alpha(x_1,...,x_n),\theta_\alpha(x_1,...,x_n))$ where $\phi_{\alpha}(\bar x) \in \mathcal L(M)$ is large and $\theta(\bar x)$ is an $\mathcal L_<$-quantifier-free formula with parameters in $M$. We build $f$ by induction on $\alpha<\lambda$ so $f$ will be defined as an increasing union $f=\cup_{\alpha<\lambda} f_\alpha$, where each $f_\alpha$ has domain a subset of $M$ of size $|\alpha|$.

At limit stages, set $f_\alpha=\cup_{\beta<\alpha} f_\beta$.

Assume $f_\alpha$ has been defined. Let $A$ be the set of parameters of $\theta_\alpha$. If needed, start by increasing $f_\alpha$ to an injection $f'$ defined on $Dom(f_\alpha) \cup A$ by defining it on $A$ in an arbitrary way. Let $\theta'$ be the formula built from $\theta_\alpha$ changing every parameter by its image by $f'$. If $\theta'$ is not large (in particular if it is inconsistant), we let $f_{\alpha+1}=f'$.

Otherwise, we can find a tuple $\bar a=(a_0,...,a_{n-1})$ in $M$ such that $M \models \phi(\bar a)$ and no $a_k$ is in the domain of $f'$. Pick any $\bar c=(c_0,...,c_{n-1})$ in $\Omega$, such that $\Omega \models \psi(\bar c)$ and no $c_k$ is in the range of $f'$. We define $f_{\alpha+1}$ as $f'\cup \{(a_k,c_k) : k<n\}$.

Once $f$ is defined, we expand $M$ to $\mathcal L'$ by letting $a<b$ if and only if $f(a) <_{\Omega} f(b)$. By construction the resulting structure is a model of $\tilde T$.
\end{proof}

\section{Counterexamples}

\subsection{$\mathbb F_p$ vector space}

We prove the following result.

\begin{thm}\label{countervspace}
If $T$ defines an infinite dimensional vector space over some prime field $\mathbb F_p$, and $T$ is linearly ordered by $<$, then $T$ has $IP$.
\end{thm}

First we need to recall the following Ramsey-type result. See for example \cite{Ramsey}, Section 2.4, Theorem 9.

Here $F$ is a fixed finite field.

\begin{thm}\label{ramseyvs}
For all $r,t,k \geq 1$, there is some $n$ such that if the $t$-dimensional affine subspaces of $F^n$ are $r$-colored, there exists a $k$-dimensional affine subspace all of whose $t$-dimensional affine subspaces have the same color.
\end{thm}

Now we prove Theorem \ref{countervspace}.

\begin{proof}
We identify $F=\mathbb F_p$ with the set $\{0,1,...,p-1\}$ and order $F$ by setting $0<1<...<p-1$.

Without loss, $T$ is just an ordered vector space over $F$. Let $M\models T$ and pick an infinite free family $\langle a_i : i<\omega\rangle$. Let $k$ be any integer and set $t=1, r=p!$. Finally, let $n$ be given by Theorem \ref{ramseyvs} for those values of $r,t,k$. We consider the vector space $A$ spanned by $(a_1,....,a_n)$ and identify it with $F^n$ using $(a_1,...,a_n)$ as canonical base. Let $\varolessthan$ be the lexicographic order on $A$ (identified with $F^n$).

Let $L\subset A$ be an affine line in $A$, $L=\{d_0,...,d_{p-1}\}$ with $d_0 \varolessthan d_1 \varolessthan ... \varolessthan d_{p-1}$. We assign to $L$ a `color' $c(L)$ from the set $S$ of permutations of $F$ in the following way: $c(L)$ is the unique permutation $\pi$ such that $$ d_{\pi(0)}< d_{\pi(1)} < .... < d_{\pi(p-1)},$$ where $<$ is the given order on $M$.

By Theorem \ref{ramseyvs}, we can find a $k$-dimensional affine subspace $W$ of $A$ all of whose $1$-dimensional affine subspaces have the same color $\pi$. Let $(b^1,...,b^k)$ be a basis of $Vect(W)$ such that for each $l\leq k$, $b^l$ is $\varolessthan$-minimal in $Vect(W)\setminus Vect(b^1,..,b^{l-1})$. For each $l\leq k$, let $\omega_l$ be the least $i$ such that the $i$'th coordinate of $b^l$ (in the basis $(a_1,...,a_n)$) is non zero. We have $$\omega_1 > \omega_2 > ... > \omega_k$$ and for each choice of $\mathfrak s=(s_1,...,s_k)\in F^k$, there is a unique element $d_{\mathfrak s}$ of $W$ such that for every $i$, the $\omega_i$'th coordinate of $d_{\mathfrak s}$ is equal to $s_i$. For $I\subseteq \{1,...,k\}$, let $d_I$ be $d_{\mathfrak s}$ for $\mathfrak s=(\varepsilon_{1,I},....,\varepsilon_{k,I})$ where $\varepsilon_{i,I}=1$ if $i\in I$ and 0 otherwise.

We let $\phi_\pi(x,y)$ be the formula:
$$\phi_\pi(x,y) = \bigwedge_{i=0}^{p-2}  x+ \pi(i). y < x + \pi(i+1). y.$$

Let $I\subseteq \{1,...,k\}$ and $l\leq k$. Consider the line $L=d_I + Vect(b^l)$. Enumerating $L$ in $\varolessthan$ increasing order gives
$$d_I - \varepsilon_{l,I} b^l \varolessthan d_I - \varepsilon_{l,I} b^l + b^l \varolessthan ... \varolessthan d_I - \varepsilon_{l,I} b^l + (p-1). b^l.$$

As $L$ has color $\pi$ we see that $\phi_{\pi}(d_I,b^l)$ holds if and only if $\varepsilon_{l,I} = 0$. This proves that the formula $\phi_\pi(x,y)$ has independence rank at least $k$. As $k$ was arbitrary, and there are only finitely many possibilities for $\pi$, there is at least one value of $\pi$ for which $\phi_\pi(x,y)$ has $IP$.

\end{proof}

\begin{rem}
In the case $p=2$, the same proof works if instead of assuming that $<$ defines a linear order, we only assume that it defines a tournament ({\it i.e.}, for all $x\neq y\in M$, exactly one of $x<y$ and $y<x$ holds).
\end{rem}

\subsection{Interpreting arithmetic}\label{bbarith}

By bounded arithmetic, we mean the (incomplete) theory $T_{arith}$ consisting of formulas true in almost all structures $(\{0,...,n\};+,\times)$. A model interprets bounded arithmetic if it has an elementary extension that interprets a model of $T_{arith}$.

\begin{thm}\label{countergraph}
There is an $\omega$-stable theory $T$, NDOP of depth 2, in a language $\mathcal L$ such that: for every model $M\models T$ and every expansion $M'$ of $M$ to $\mathcal L'=\mathcal L\cup \{<\}$ in which $<$ defines a linear order, $M'$ interprets bounded arithmetic.
\end{thm}

We take as language $\mathcal L=\{E,S,R\}$ where $E$ and $S$ are binary predicates and $R$ is quaternary. The predicate $E$ defines an equivalence relation on the structure and each $E$-class is a made by $(S, R)$ into a regular colored graph as explained later.

We will define $T$, by constructing its prime model $M_0$ built as a disjoint union of finite regular graphs. Each of those graphs will be exactly one $E$-class. The $n$'th graph will contain no cycle of length $\leq n$. We will choose the finite graphs in such a way that no mater what order is put on them, the $n$'th graph interprets $(\{1,...,n\},+,\times)$, by an explicit formula not depending on $n$. The condition about cycles will ensure that the limit theory of those graphs is $\omega$-stable. 

The first observation is that, for every $n$, there is a structure $(\{1,..,N\},<_1,<_2)$ where both $<_1$ and $<_2$ are linear orders (call this a \emph{bi-order}) that interprets $(\{1,...,n\},+,\times)$. Furthermore, the formulas involved in the interpretation do not depend on $n$. This is an easy exercise whose solution is given in the appendix. As a consequence, the problem is reduced to that of interpreting two linear orders.
\\

We will work with colored graph. For us a colored graph is a structure $(G; S, R)$ where $S(x,y)$ is symmetric anti-reflexive and $R(x_0,y_0,x_1,y_1)$ is a quaternary relation which defines an equivalence relation on pairs $(x,y)\in S$. It should be thought of as saying that $\{x_0,y_0\}$ and $\{x_1,y_1\}$ are two edges in the graph of the same color. We will consider only regular colored graph, namely such that each vertex is part of exactly one edge of each color. For simplicity, we will introduce the imaginary sort $C$ of colors defined by the quotient of $\{(x,y)\in G^2 : \models S(x,y)\}$ by the equivalence relation $(x_0,y_0) E (x_1,y_1) \iff R(x_0,y_0,x_1,y_1)$. If $k\in C$ is a color and $x\in G$ is a vertex, by the $k$-neighbor of $x$, we mean the unique $y\in G$ such that $\{x,y\}$ is an edge of color $k$.

Denote by $\mathcal L_{cg}$ the language $\{S,R\}$ of colored graph. We will use $\mathcal L'_{cg}$ to denote the expanded language $\{S,R,<\}$ where $<$ is a binary predicate. Also let $\Phi(x,y,u,v)$ be the $\mathcal L'_{cg}$ formula saying that $\{u,v\}$ is an edge and if $k$ is its color, then the $k$-neighbor of $x$ is $<$-less than the $k$-neighbor of $y$.

In the rest of this subsection, we prove the following proposition.

\begin{prop}\label{randomgraph}
Let $N,c\geq 3$ be integers, let also $<_1,<_2$ be two linear orders on $N=\{0,...,N-1\}$. Then for every even integer $n$ big enough, there is a finite regular colored graph $G_n$ on $n$ vertices such that for every expansion of $G_n$ to $\mathcal L'_{cg}$ making $<$ into a linear order, there are $a,b,u,v \in G_n$ such that the structure $([a,b];<,\Phi(x,y,u,v))$ is isomorphic to $(N;<_1,<_2)$. Furthermore, $G_n$ has no cycle of length $\leq c$.
\end{prop}

Here $[a,b]$ denotes the interval with end points $a$ and $b$ in the sense of $<$.
\\

We will build the graph $G_n$ by a random procedure and show that with positive probability, we obtain what we want. Actually, we will start by building a colored multigraph $G(\sigma)$ and then modifiy it to make into an actual regular graph $G'(\sigma)$. By a colored multigraph, we mean a colored graph in which there can be two or more edges (of different colors) between two given vertices.

Let $n$ be even and big enough (we will see during the construction what big enough means). Fix some $1-\frac 1{3c} < \alpha <1$ and let $d=n^{1-\alpha}$. We take $n=\{0,...,n-1\}$ as set of vertices and $C=\{0,...,d-1\}$ as set of colors. So our final graph will be $d$-regular.

Let $\mathfrak S_n$ denote the symmetric group on $n$ elements. Our space of events is $\Omega = \mathfrak S_n^d$ equipped with the uniform probability law. Let $\sigma=(\sigma_i : i<d)$ be an element of $\Omega$. We define a colored multigraph $G(\sigma)$ as follows: The vertex set of $G(\sigma)$ is $\{0,...,n-1\}$. The set of colors is $\{0,...,d-1\}$. For $k<d$, we draw an edge of color $k$ between vertices $a$ and $b$ if and only if for some $l< n/2$, we have $\{\sigma_k(2l),\sigma_k(2l+1)\}=\{a,b\}$.

A cycle of length $r\geq 2$ is a sequence $(a_0,...,a_{r-1})$ of distinct vertices and a sequence $(e_0,...,e_{r-1})$ of distinct edges such that for each $i<r$, $e_i$ is an edge between $a_i$ and $a_{i+1}$ (addition is modulo $r$). In particular, a cycle of length 2 in $G(\sigma)$ is given by two vertices and two different edges linking them. A cycle is said to be \emph{small} if it is of length $\leq c$.

If $k\in C$ is a color, we define $\lessdot_k$ by $x \lessdot_k y$ if and only if the $k$-neighbor of $x$ is $<$-less than the $k$-neighbor of $y$. We will show that each of the following events occurs with probability converging to 1 as $n$ tends to $+\infty$:

\begin{enumerate}
\item The number of small cycles in $G(\sigma)$ is less than $d^{2c}$, and we can obtain a regular graph $G'(\sigma)$ with no small cycles by changing at most $2.d^{2c}$ edges of $G(\sigma)$,

\item For every ordering $<$ of the vertices of $G(\sigma)$, we can find at least $d^{3c}$ values of $(a,b,k)\in G^2 \times C$ such that $([a,b];<,\lessdot_k)$ is isomorphic to $(N;<_1,<_2)$.
\end{enumerate}

\subsubsection{Removing small cycles}

We perform some surgery to remove small cycles from the multigraph $G(\sigma)$ and obtain a regular graph $G'(\sigma)$ with no small cycles.

We compute the expectancy of the number of small cycles in $G(\sigma)$. Let $2\leq s\leq c$. Let $H=(V_H,E_H)$ be the graph consisting of a unique cycle of length $s$. Fix some $f:V_H \rightarrow \{0,...,n-1\}$ injective and some $g:E_H \rightarrow \{0,...,d-1\}$. Assume that for any two edges $e,e'\in E_H$ having a vertex in common, $g(e) \neq g(e')$.
Call such a pair $(f,g)$ a subgraph of $G(\sigma)$ if for every edge $e=\{a,b\}\in E_H$, $\{f(a),f(b)\}$ is an edge of $G(\sigma)$ and has color $g(e)$. For a given edge $e$, the probability $p_0$ that this occurs is $$p_0 = \frac 1{n-1}.$$

If $g$ is injective, then all those events are independent, so the probability $p$ that $(f,g)$ is a subgraph satisfies $$p \sim \left (\frac 1 {n} \right )^s.$$
If $g$ is not injective, the events are not independent, but nevertheless, looking color by color, the same estimate can easily be seen to be true.

The number of such pairs $(f,g)$ is less than $n^s d^s$. So the expected value of the number of cycles of length $s$ is asymptotically at most $$n^s d^s \left (\frac 1 {n}\right )^s= d^s.$$

Let $X$ denote the number of small cycles in $G(\sigma)$. Summing over all $s\leq c$, we see that $E(X)= O(1).d^c$. By the first moment method $$Prob(X\geq d^{2c}) \leq O(1).d^{-c}.$$ In particular, with probability converging to 1 at $n$ tends to $+\infty$, the number of small cycles in $G(\sigma)$ is less than $d^{2c}$. Consider a $\sigma$ that has that property.

We now modify the graph $G(\sigma)$ so as to remove all small cycles. We show that this can always be done by changing at most $2.d^{2c}$ edges.

Consider each small cycle one after the other. Let $C$ be such a cycle, of size $s\leq c$. Pick any edge $\{a,b\}$ in $C$. Assume it is drawn with color $k$. We will choose an edge $\{a',b'\}$ drawn with color $k$, erase both those edges, and draw instead, with color $k$, the edges $\{a,a'\}$ and $\{b,b'\}$. We have to choose $\{a',b'\}$ so that no new small cycle is introduced. It is easy to check that this will happen if the following two properties are satisfied:

-- the graph distance between $a$ and $a'$ is at least $c+2$,

-- the edge $\{a',b'\}$ does not belong to a small cycle.

Note that the first condition implies that the distance between $b$ and $b'$ is at least $c$. We see that the number of edges $\{a',b'\}$ which fail to satisfy those properties is at most $O(1).(d^{c+2}+d^{2c})<n$, so we can find $a',b'$ as required.

In the end, we have modified at most $2d^{2c}$ edges. Call $G'(\sigma)$ the graph obtained at the end of the surgery.

\subsubsection{Interpreting bi-orders}

All is left to prove is that if $(N;<_1,<_2)$ is given in advance, we can interpret it with high probability in $G(\sigma)$ (uniformly in $N$) in many different ways. This will automatically imply that we can also interpret it also in $G'(\sigma)$.

Recall that the binomial distribution with parameters $n,p$ is defined as the distribution of the sum of $n$ independent random variables being equal to 1 with probability $p$ and to 0 with probability $1-p$. We will need the following fact about binomial distributions. It is a special case of Hoeffding's inequality (or of Chernoff's bound). See for example \cite{Hoeffding}.

\begin{lemme}\label{binom}
Let $X$ be a random variable whose law is a binomial distribution with parameters $n,p$. Then for any $x$ we have:

$$Prob(X\leq x) \leq \exp \left ( -2 \frac {(np-x)^2}{n} \right ).$$
\end{lemme}

We are given a fixed structure $(N;<_1,<_2)$. Let $<$ be any order on $\{0,...,n-1\}$ (the set of vertices of $G(\sigma)$), without loss, the usual order. Given $a<n$ and $k\in C$ a color, define the event $\mathcal E_{a,k}$ as $$([a,a+N-1]; <, \lessdot_k) \text{ is isomorphic to }(N;<_1,<_2).$$

Fix a value of $k$. Let $D\subset \{0,...,n-1\}$ be a subset of size $<n/4N$ closed under taking the $k$-neighbor. Assume we know the values of $\sigma_k(x)$ for every $x\in D$ and the values of $\sigma_l(x)$ for $l<k$ and every $x$. Take some $a<n$ such that $D\cap \{a,a+1,...,a+N-1\}=\emptyset$. We look at the values of  the $k$-neighbors of the points $\{a,a+1,...,a+N-1\}$. With probability $>1/2$ those neighbors are disjoint from $\{a,a+1,...,a+N-1\}$ so we are considering $2N$ distinct points. Consider the bi-order $([a,a+N-1];<,\lessdot_k)$. Clearly, all $N!$ bi-order have the same probability of occurring. So, knowing $\sigma_k$ restricted to $D$ and $\sigma_l$ for $l<k$, the event $\mathcal E_{a,k}$ holds with probability at least $\frac 1{2N!}$.

Assume we know the values of $\sigma_l$ for $l<k$. We start with $D_0=\emptyset$, take $a_0=0$. With probability at least $\frac 1{2N!}$, the event $\mathcal E_{a_0,k}$ holds. Next we let $D_1$ contain $\{a_0,a_0+1,...,a_0+N-1\}$ and all $k$-neighbors of those points. We take some $a_1$ as in the previous paragraph. Again with probability at least $\frac 1{2N!}$, $\mathcal E_{a_1,k}$ holds. We can iterate this at least $\frac n{8N^2}$ times (as $D$ increases by at most $2N$ each time). We do this for each color, one color after the other. Let $Y$ be the number of events $\mathcal E_{a,k}$ that hold. Then for every $x$, $$Prob(Y \leq x) \leq Prob(B \leq x)$$ where $B$ is a random variable whose law is a binomial distribution with parameters $\frac {dn}{8N^2}$, $\frac 1{2N!}$.

By \ref{binom}, the probability that at least $d^{3c}$ of the events $\mathcal E_{a,k}$ succeed is at least $$1-\exp \left(-16N^2 \frac{(\frac {dn}{16N^2N!} - d^{3c})^2}{nd}\right ) = 1-\exp \left ( -n^{2-\alpha}.O(1)  \right ),$$

As there are $n!\leq \exp(n \ln n)$ different orders on $\{0,...,n-1\}$, we see that with probability converging to 1 as $n$ tends to $+\infty$, for every ordering $<$ of the vertices, there are at least $d^{3c}$ values of $(a,k)\in n \times C$ for which $$([a,a+N-1];<,\lessdot_k)\text{ is isomorphic to }(N;<_1,<_2).$$

In particular, if furthermore $G(\sigma)$ has less than $d^{2c}$ small cycles, there is a choice of $(a,k)$ such that no edge of color $k$ having an end point in $\{a,...,a+N-1\}$ is changed during the construction of $G'(\sigma)$. Therefore with this choice of $(a,k)$, we obtain an interpretation of $(N;<_1,<_2)$ in $G'(\sigma)$ as $$([a,a+N-1],<,\Phi(x,y,u,v))$$ where $\{u,v\}$ is any edge of color $k$.

This ends the proof of \ref{randomgraph}.

\subsection{The full structure}

We now conclude with the proof of Theorem \ref{countergraph}.

We recall that we set $\mathcal L=\{E,S,R\}$ where $E$ is a binary relation. Fix an $3 \leq n<\omega$. Let $(N;<_1,<_2)$ be the structure $P_n$ given by \ref{arith}. Let also $(G_n;S,R)$ be the colored graph given by Proposition \ref{randomgraph} with $c=n$. We define the $\mathcal L$-structure $M_0$ as follows. The reduct to $\{S,R\}$ is just the disjoint union of the graphs $(G_n;S,R)$. The predicate $E$ is interpreted as an equivalence relation, such that two points $x,y \in M_0$ are $E$-equivalent if and only if they come from the same $G_n$.

What does $Th(M_0)$ look like? Let $M$ be elementary equivalent to $M_0$. Then in $M$, $E$ defines an equivalence relation which has $\omega$ finite classes. The substructure of $M$ formed by the union of those finite classes is isomorphic to $M_0$. An infinite class of $M$, equipped with $(S,R)$ is a regular colored graph with no cycles and infinite degree. So as a graph, it is a union of trees. Finally, there are no $S$-relations between points in different classes.

We see that to every $E$ class are associated two regular types: one for a new connected component and one for a new color. The type of a new $E$ class is itself a regular type. A model is entirely determined up to isomorphism by the number of infinite $E$-classes and for each such class, the number of connected components and the number of colors. So easily, the resulting theory is $\omega$-stable, NDOP of depth 2.

Finally, if $<$ is an expansion of $M_0$ by a linear order, then in the $k$'th equivalence class, we can interpret $(\{0,...,k\};+,\times)$, uniformly in $k$. In particular the expansion interprets bounded arithmetic.

 \section*{Appendix: interpreting arithmetic from two linear orders}

The following is certainly well known, but we include it for completeness.
\begin{prop}\label{arith}
There are formulas $\phi(x,y,z;\bar t)$, $\psi(x,y,z;\bar t)$ in the language $\mathcal L=\{<_1,<_2\}$ such that for each $k<\omega$, there is an $\mathcal L$-structure $P_k$ satisfying :
\begin{itemize}
\item $<_1$ and $<_2$ define linear orders on $P_k$,
\item there are $a_0,a_1,\bar a$ in $P_k$ such that the structure $$([a_0,a_1]_2;\phi(x,y,z;\bar a),\psi(x,y,z;\bar a))$$ is isomorphic to $(\{0,...,k-1\},+,\times)$.\\(where $[a,b]_2$ denotes the interval $a \leq_2 x \leq _2 b$).
\end{itemize}
\end{prop}
\begin{proof}
Let $k<\omega$ and set $n=10.k^2$. Let $P_k$ have universe $\{0,...,n-1\}$ and let $<_1$ be the usual order on that set. For $i\in \{0,...,k-1\}$, let $b_i = 7k.(i+1)$. Let $a_1=b_{k-1}$. The points of $P_k$ greater then $b_{k-1}$ will be called \emph{delimiters}.

Now we construct $<_2$ in the following way : an initial segment is $b_0 <_2 b_1 <_2 ... <_2 b_{k-1}$. Those points will correspond to $\{0,..,k-1\}$ in the interpretation. For $l<n$, $r<k$ we say that $l$ codes for $r$ if $b_{r-1} < l < b_r$ (setting $b_{-1}=-1$). Now we encode the graph of addition. Let $(A_i)_{i<N}$ be an enumeration of all triples $(r,s,t)\in \{0,..,k-1\}^3$ such that $r+s=t$. Build a map $f: N \rightarrow P_k^4$ such that for each $i<N$, if $f(i)=(c_r,c_s,c_t,d)$ then $c_r,c_s$ and $c_t$ code for $r,s,t$ respectively and $d=a_0+i+1$ is a delimiter. Also impose that no point appears in two tuples $f(i),f(i')$ for $i\neq i'$. Now set $<_2$ such that $c_r <_2 c_s <_2 c_t <_2 d$ each time $f(i)=(c_r,c_s,c_t,d)$, and those four points form a interval of $<_2$. Place those intervals one after another in any order. Let $a_2$ be the last delimiter placed. Next do the same for the graph of multiplication, and let $a_3$ be the last delimiter placed there. Finally place the elements that are left in any order after $a_3$.

Now it is easy to check that one can interpret $(\{0,...,k-1\};+,\times)$ in this structure, with parameters $a_1,a_2,a_3$, and that the formulas involved do not depend either on the choices made nor on $k$.

\end{proof}

\end{document}